\documentclass[12pt]{article}
\usepackage{latexsym}
\usepackage{tikz}
\usepackage{tkz-graph}
\usetikzlibrary{shapes}

\newtheorem{Theorem}{Theorem}[section]

\newtheorem{Lemma}[Theorem]{Lemma}
\newtheorem{Observation}[Theorem]{Observation}

\def\inst#1{$^{#1}$}
\def\PP{($P_{5},\overline{P}_5$)}

\title{On color-critical ($P_{5},\overline{P}_5$)-free graphs}

\author{
  Harjinder S. Dhaliwal\inst{1}
  \and Ang\`ele M. Hamel\inst{1}
  \and Ch\'inh T. Ho\`ang\inst{1}
  \and Fr\'ed\'eric Maffray\inst{2}
  \and Tyler J. D. McConnell\inst{1}
  \and Stefan A. Panait\inst{1}
}

\begin{document}
\maketitle

\begin{center}
{\footnotesize

\inst{1} Department of Physics and Computer Science, Wilfrid Laurier
University, \\ Waterloo, Ontario, Canada

\inst{2} CNRS, Laboratoire G-SCOP, UJF, INPG, Grenoble, France}

\end{center}
\begin{abstract}
A graph is $k$-critical if it is $k$-chromatic but each of its proper
induced subgraphs is ($k-1$)-colorable.  It is known that the number
of $4$-critical $P_5$-free graphs is finite, but there is an infinite
number of $k$-critical $P_5$-free graphs for each $k \geq 5$.  We show
that the number of $k$-critical $(P_5, \overline{P}_5)$-free graphs is
finite for every fixed $k$.  Our result implies the existence of a
certifying algorithm for $k$-coloring $(P_5, \overline{P}_5)$-free
graphs.

\end{abstract}

\section{Introduction}\label{sec:intro}

Graph coloring is a well-studied problem in computer science and
discrete mathematics.  Determining the chromatic number of a graph
is a NP-hard problem.  But for many classes of graphs, such as
perfect graphs, the problem can be solved in polynomial time.
Recently, much research has been done on coloring $P_5$-free
graphs.  Finding the chromatic number of a $P_5$-free graphs is
NP-hard \cite{KraKra2001}, but for every fixed $k$, the problem of
coloring a graph with $k$ colors admits a polynomial-time
algorithm \cite{HoaKam2008, HoaKam2010}.  Research has also been
done on ($P_5, \overline{P}_5$)-free graphs (graphs without $P_5$
and its complement $\overline{P}_5$).  In \cite{GiaRus1997}, a
polynomial-time algorithm is found for finding an approximate
weighted coloring of a ($P_5, \overline{P}_5$)-free graph
(definitions not given here will be given later).  Recently,
\cite{HoaLaz2013} gave a polynomial time algorithm for finding a
minimum weighted coloring of a ($P_5, \overline{P}_5$)-free graph.

The algorithms in \cite{HoaKam2008, HoaKam2010, HoaLaz2013} produce a
$k$-coloring if the input graph is $k$-colorable.  However, when the
graph is not, the algorithms do not produce an easily verified
certificate for a ``NO'' answer.  The point of view in this article is
motivated by the idea of a ``certifying algorithm''.  An algorithm is
{\it certifying} if it returns with each output a simple and easily
verifiable certificate that the particular output is correct.  For
example, a certifying algorithm for the bipartite graph recognition
would return either a 2-coloring of the input graph, thus proving that
it is bipartite, or an odd cycle, thus proving it is not bipartite.  A
certifying algorithm for planarity would return either an embedding of
the graph in a plane, or one of the two Kuratowski subgraphs proving
the input graph is not planar.

A graph is $k$-critical if it is $k$-chromatic but each of its proper
induced subgraphs is $(k-1)$-colorable.  In \cite{BruHoa2009} and also
\cite{MafMor}, a certifying algorithm for $3$-colorability of
$P_5$-free graphs is provided by showing that the number of
$4$-critical $P_5$-free graphs is finite.  Given this result, one may
ask whether the same statement holds for $k$-critical $P_5$-free
graphs for any $k \geq 4$.  However, \cite{HoaMoo2013} shows the
number of $k$-critical $P_5$-free graphs is infinite for $k \geq 5$.
Here we prove that the number of $k$-critical ($P_5,
\overline{P}_5$)-free graphs is finite for every fixed $k$, thereby
establishing a certifying algorithm for $k$-colorability of ($P_5,
\overline{P}_5$)-free graphs.  In section~\ref{sec:definitions}, we
give definitions and background on our problem.  In
section~\ref{sec:characterization}, we give the proofs of our main
results.  In section~\ref{sec:5-critical}, we give the exact number of
$k$-critical {\PP}-free graphs for $k \leq 8$; in particular, we will
construct a list of all 5-critical {\PP}-free graphs.

\section{Definitions and background}\label{sec:definitions}
A $k$-coloring of a graph $G=(V,E)$ is a mapping $f: V \rightarrow
\{1,\ldots, k\}$ such that $f(u) \not= f(v)$ whenever $uv \in E$.
Given a coloring, a {\it color class} is the set of all vertices
of the same color.  The chromatic number $\chi(G)$ of a graph $G$
is the smallest $k$ such that $G$ is $k$-colorable.  $G$ is
$k$-chromatic if $\chi(G) = k$.  A graph $G$ is {\em $k$-critical}
if it is $k$-chromatic and none of its proper induced subgraphs is
$k$-chromatic (that is, all of its proper induced subgraphs are
$(k-1)$-colorable).  We say that a graph is {\it critical} if it
is $k$-critical for some $k$.  Let $N(v)$ be the set of neighbors
of $v$. A set $X$ of vertices of a graph $G=(V,E)$ is a {\em
module} if for all $v \not\in X$, $X \subseteq N(v)$ or $N(v) \cap
X = \phi$.  Module $X$ is {\em trivial} if $|X| = 1$ or $X = |V|$.
Unless otherwise stated, a module in this paper is non-trivial.  A
vertex of $G$ is {\it universal} if it is adjacent to every other
vertex of $G$. Vertices $u,v$ are comparable if $N(u) \subseteq
N(v)$, or vice versa. If $X$ is a set of vertices of $G$, then
$G[X]$ denotes the subgraph if $G$ induced by $X$.  A set $A$ of
vertices is {\it complete} to a set $B$ of vertices if there are
all edges between $A$ and $B$.  Given two graphs $G$ and $ H$, the
graph $F$ is the {\it join} of $G$ and $H$ if $F$ is obtained by
taking $G$ and $H$ and joining every vertex in $G$ to every vertex
in $H$ by an edge.  As usual, $K_t$ denotes the clique on $t$
vertices; and $C_t$ denotes the induced cycle on $t$ vertices. The
complement of $G$ is denoted by $\overline{G}$.

Let $G$ be a graph with a module $M$, where $M$ needs not be
non-trivial.  Consider the graph $G'$ obtained from $G$ by
removing $M$, adding another graph $H$, and adding all edges
between a vertex $x \in G-M$ and all vertices of $H$ whenever $x$
has a neighbor in $M$ in the graph $G$.  We say that $G'$ is
obtained from $G$ by {\it substituting} $H$ for $M$.  A {\em buoy}
is a graph $B$ obtained from a $C_5$ by, for each vertex $x$ of
the $C_5$, substituting a graph $B_x$ for $x$.  The graph $B_x$ is
a {\em bag} of the buoy.  Thus, the buoy will have $5$ bags, and
we label these $B_1, B_2, B_3, B_4, B_5$ in the cyclic order.
Without loss of generality, we will start with $B_1$ when
discussing bags.  Note each $B_i$ is a module of $B$.  The buoy
$B$ is {\it full} if $B$ contains all vertices of $G$; $B$ is {\it
universal} if every vertex in $G-B$ is universal; $B$ is a {\it
join buoy} if every vertex of $B$ is adjacent to every vertex of
$G-B$.  A \emph{pseudo-buoy} is defined as a buoy except that any
of the sets $B_1, \ldots, B_5$ may be empty.  The following two
lemmas are folklore and are easy to establish.

\begin{Lemma}\label{lem:neighbor}
In a critical graph there do not exist two comparable vertices.
$\Box$
\end{Lemma}

\begin{Lemma}\label{lem:connected}
A critical graph is connected. $\Box$
\end{Lemma}
The following lemma is also well known and we will rely on it for
our proofs.
\begin{Lemma}\label{lem:cliqueCutset}
A critical graph that is not a clique does not contain a clique
cutset.
\end{Lemma}
{\it Proof.} Let $G$ be a $k$-critical graph with a clique cutset
$C$.  Thus, $G-C$ can be partitioned into two sets $A, B$ such
that there are no edges between $A$ and $B$.  Let $G_A$ (resp.,
$G_B$) be the subgraph of $G$ induced by $C \cup A$ (resp., $C
\cup B$).  Consider a coloring of $G_A$ (resp., $G_B$) with
$\chi(G_A)$ (resp., $\chi(G_B)$) colors.  Since $G$ is critical,
we have $\chi(G_A) < k$ and $\chi(G_B) < k$.  Without loss of
generality, we may assume $\chi(G_A) \geq \chi(G_B)$.  Then a
$\chi(G_A)$-coloring of $G$ can be obtained by identifying the
colors of $G_A$ and $G_B$ on $C$.  But this implies $\chi(G) =
\chi(G_A) < k$, a contradiction. $\Box$

\begin{Lemma}\label{thm:FG}{\rm \cite{FouGia1995}}
Given a $(P_5, \overline{P}_5)$--free graph $G$, every $C_5$ of $G$ is
contained in a buoy which is either full or a module of $G$.  $\Box$
\end{Lemma}

Let $\omega(G)$ denotes the {\it clique number} of $G$, ie. the number of vertices in a largest clique of
$G$. A graph $G$ is {\it perfect} if for each induced subgraph $H$ of $G$,
we have $\chi(H)=\omega(H)$.
\begin{Observation}\label{obs:perfect}
If $G$ is perfect and $k$-critical,
then $G$ is the graph $K_k$, the clique on $k$ vertices.
\end{Observation}
{\it Proof.} Let $G$ be a perfect $k$-critical graph. Thus, we have $\chi(G)=\omega(G)=k$.
For a vertex $v$ of $G$, we have
$\chi(G-v) = k -1$ by definition. Since $G$ is perfect, we have $\omega(G-v)=k-1$. It follows that any largest clique of $G$ must contain $v$. Since $v$ was chosen arbitrarily, a largest clique of $G$ must contain all vertices of $G$, that is, $G$ is a clique on $k$ vertices. $\Box$

We will also need the following result.
\begin{Lemma}\label{lem:CHMW}{\rm \cite{ChvHoa1987}}
$(P_5, \overline{P}_5, C_5)$-free graphs are perfect.  $\Box$
\end{Lemma}

\section{The structure of $k$-critical {\PP}-free
graphs}\label{sec:characterization}
In this section, we study properties of maximal buoys in
{\PP}-free graphs and we show that every critical {\PP}-free graph
$G$ can be constructed by two simple operations: (i) $G$ is the
join of two smaller critical graphs, or (ii) $G$ is obtained from
the $C_5$ by substituting some smaller critical graphs for each
vertex on the $C_5$. As a consequence, we obtain a certifying
algorithm for $k$-colorability of {\PP}-free graphs.
\subsection{Structural results}

We begin by  establishing a number of preliminary results.
\begin{Lemma}\label{lem:module-critical}
If $G$ is $k$--critical with a non-trivial module $M$, then $M$ is
$\ell$--critical for some $\ell < k$.
\end{Lemma}
\noindent {\it Proof}.  As $M$ is a module, we can partition the
vertices of $G$ into three sets: $M$; $N$, the set of vertices in
$G-M$ adjacent to the vertices of $M$; and $R$, the set of
vertices in $G-M$ having no neighbors in $M$.

Let $\chi(M)=\ell$ for some $\ell\leq k$.  If $\ell=k$, then $N$ must
be empty since any color assigned to a vertex in $N$ must be different
from the colors assigned to vertices in $M$ (and so $G$ would require
more than $k$ colors).  This would imply that $G$ is not connected, a
contradiction to Lemma~\ref{lem:connected}.  So, we know $\ell < k$.

Now, we need to  show that $M$ is $\ell$-critical; that is, we must
show that  $M-y$ is $(\ell-1)$-colorable for any vertex $y\in M$.

Suppose on the contrary that there is a vertex $x \in M$ with $\chi(M-x)
= \ell$.  Now consider a $(k-1)$-coloring $\beta$ of the graph $G-x$
which we know must exist because $G$ is $k$-critical.  Then at least
$\ell$ colors of $\beta$ appear in $M-x$, because $\chi (M-x) = \ell$.
Let $C_M$ be the set of colors of $\beta$ that appear in $M$, and let
$C_N$ be the set of colors of $\beta$ that appear in $N$.  We have
$C_M \cap C_N = \emptyset$ since there are all edges between $M$ and
$N$.

Let $\theta$ be an $\ell$-coloring of $M$ that uses the colors of
$C_M$; such a coloring exists because $M$ is $\ell$-colorable. Now
consider an assignment $\gamma$ of colors to $G$ which uses the
same colors for $R$ and $N$ as in $\beta$ but uses the coloring of
$\theta$ for $M$, that is, for a vertex $x$, if $x \in R \cup N$,
then $\gamma(x) = \beta(x)$; and if $x \in M$, then $\gamma(x)=
\theta(x)$. Then $\gamma$ is a valid coloring because for any two
adjacent vertices $x,y$, we have $\gamma(x) \not= \gamma(y)$; in
particular, if $x \in M$ and $y \in N$ then $\gamma(x) \in C_M$
and $\gamma(y) \in C_N$ and so $\gamma(x) \not= \gamma(y)$.  So,
$\gamma$ uses the same number of colors as $\beta$; that is, $G$
is $(k-1)$-colorable, a contradiction. $\Box$

\begin{Lemma}\label{lem:join-critical}
Let $G=(V,E)$ be any graph.  Suppose that $V$ admits a partition into
two non-empty sets $V_1$ and $V_2$ such that $V_1$ is complete to
$V_2$.  Then $G$ is critical if and only if the two graphs $G[V_1]$
and $G[V_2]$ are critical.
\end{Lemma}
\noindent{\it Proof.} Let $k=\chi(G)$.  Write $G_i=G[V_i]$ for
each $i\in\{1,2\}$.

First suppose that $G$ is $k$-critical.  If $V_i$ is a single vertex
then it is 1-critical; otherwise, $V_i$ is a module of $G$ and so by
Lemma~\ref{lem:module-critical}, $G[V_i]$ is critical.

Now suppose that each of $G_1$ and $G_2$ is critical.  Let $k_1=
\chi(G_1)$ and $k_2= \chi(G_2)$.  So $k=\chi(G)= \chi(G_1)+ \chi(G_2)=
k_1+k_2$.  Pick any $x\in V_1$.  Since $G_1$ is critical, $G_1 - x$
admits a $(k_1-1)$-coloring.  We can combine this coloring with any
$k_2$-coloring of $G_2$, using a disjoint set of colors, to obtain a
$(k-1)$-coloring of $G$.  The same holds if $x\in V_2$.  So $G$ is
critical.  $\Box$

\medskip

\begin{Lemma}\label{lem:substitute-critical}
Let $G=(V,E)$ be a $k$-critical graph for some $k$.  Let $M$ be a
module of $G$ with $\chi(M) = \ell$ for some $\ell$.  Let $G'$ be the
graph obtained from $G$ by substituting a clique $K$ on $\ell$
vertices for $M$.  Then $G'$ is also $k$-critical.
\end{Lemma}
{\it Proof.} It is easy to see that $\chi(G') = \chi(G) = k$. We
need to prove that every proper induced subgraph of $G'$ is
$r$-colorable for some $r < k$.  Let $H'$ be a proper induced
subgraph of $G'$.  Let $N$ be the set of vertices of $G$ with some
neighbor in $M$ and let $R = V - N$. With respect to $G'$, the
sets $N$ and $R$ remain unchanged (vertices in $N$ would have
neighbors in $K$).

Suppose some vertex $x$ in $K$ does not belong to $H'$.  Let $t$ be
the number of vertices of $K$ that are in $H'$ with $t < \ell$.  In
$M$, consider an induced subgraph $M_t$ with chromatic number $t$.
Such graph exists since $t < \ell = \chi(M)$.  Consider the subgraph
of $G$ induced by $N \cup R \cup M_t$.  It admits an $r$-coloring
$\alpha$ for some $r < k$.  In this coloring, at least $t$ colors
appear in $M$.  From $\alpha$, we can construct an $r$-coloring of
$H'$ by (i) for $v \in H' \cap (N \cup R)$, giving $v$ the color
$\alpha(v)$, (ii) giving each of the $t$ vertices of $K \cap H'$ a
distinct color used by $\alpha$ on $M$.

Thus, we may assume all vertices of $K$ belong to $H'$.  It follows
that some vertex $x \in N \cup R$ is not in $H'$.  Let $H$ be the
proper induced subgraph of $G$ obtained from $H'$ by substituting $M$
for $K$.  Since $H$ admits an $r$-coloring $\alpha$ for some $r < k$,
we may obtain an $r$-coloring of $H'$ from $\alpha$ by (i) for $v \in
H' \cap (N \cup R)$, giving $v$ the color $\alpha(v)$, (ii) giving
each of the vertices of $K$ a distinct color that $\alpha$ uses on
$M$.  $\Box$

\begin{Observation}\label{obs:substitution}
Let $G$ and $H$ be two  {\PP}-free graphs. Suppose $G$ contains  a
module $M$. Let $G'$ be the graph obtained from $G$ by
substituting $H$ for $M$. Then $G'$ is {\PP}-free.
\end{Observation}

\noindent {\it Proof.} Let $N$ be the set of vertices of $G-M$
with some neighbours in $M$, and let $R=V(G)-M-N$. Then $V(G') = R
\cup N \cup V(H)$ and $H$ is a module of $G'$. Suppose $G'$
contains a graph $P$ which is isomorphic to a $P_5$ or
$\overline{P}_5$. Observe that $P$ contains no module. For any
vertex $x \in M$ and any vertex $y \in H$, the subgraph $G_x$ of
$G$ induced by $R \cup N \cup \{x\}$ is isomorphic to the subgraph
$G_y'$ of $G'$ induced by $R \cup N \cup \{y\}$. We have $ \left
\vert V(P) \cap H \right \vert \geq 2$, for otherwise $P$ lies
completely in $G_x$ and therefore in $G$, a contradiction. We also
have $V(P) - H \not= \emptyset$, for otherwise $P$ lies entirely
in $H$, a contradiction. But then in $P$, the set $V(P) \cap H$ is
a module, a contradiction. $\Box$

\begin{Lemma}\label{lem:buoy}
Let $G$ be a critical {\PP}-free graph.  If $G$ contains a $C_5$, then  $G$ contains a
full or a join buoy.
\end{Lemma}

\noindent {\it Proof.} The smallest graph with a buoy is the $C_5$
for which the Lemma obviously holds. We now prove the Lemma by
contradiction. Let $G=(V,E)$ be a smallest $k$-critical {\PP}-free
graph without a join or full buoy. Let $B$ be a maximal buoy of
$G$. (Here, ``Maximal'' is meant with respect to set-inclusion. In
particular, a maximal buoy may not be a largest buoy.) By
Theorem~\ref{thm:FG}, we know $B$ is a module of $G$. Let $N$ be
the set of vertices of $G-B$ with some neighbors in $B$ and let $R
= V - N - B$.  We know $R \not= \emptyset$, for otherwise $B$ is a
join buoy, a contradiction.  Let $G'$ be the graph obtained from
$G$ by substituting a clique $K$ on $\chi(B)$ vertices for $B$. By
Lemma~\ref{lem:substitute-critical}, we know $G'$ is also
$k$-critical. Note that $K$ is a module of $G'$, and with respect
to $K$ and $G'$, the sets $N$ and $R$ remain unchanged (vertices
in $N$ would have neighbors in $K$). By
Observation~\ref{obs:substitution}, $G'$ is {\PP}-free. The graph
$G'$ is not a clique because there is a non-edge between $R$ and
$K$. By Observation~\ref{obs:perfect}, $G'$ is not perfect. Now,
Lemma~\ref{lem:CHMW} implies $G'$ contains a $C_5$. Note that $G'$
has fewer vertices than $G$ ($B$ has a non-edge but $K$ is a
clique).  The minimality of $G$ implies $G'$ contains a buoy $B'$
that is full or join.

We now prove $K-B' \not= \emptyset$. Suppose $K$ is completely
contained in $B'$. Since $K$ is a module of $G'$, $K$ lies
entirely in one bag of $B'$.  Let $B''$ be the  graph obtained
from $B'$ by substituting $B$ for $K$. Then $B''$ is a buoy of $G$
and strictly contains $B$, a contradiction to the maximality of
$B$. So, $K$ must contain some vertex of $G'-B'$. This implies
$B'$ is not a full buoy of $G'$.

Now, we may assume $B'$ is a join buoy of $G'$, that is, $G'$ can
be partitioned into two set $F_1 = B'$ and $F_2$ such that there
are all edges between $F_1$ and $F_2$.  The clique $K$ cannot
contain a vertex $b$ in $B'$ and and a vertex $a$ in $F_2$ because
some vertex of $B'$ would be non-adjacent to  $b$ and adjacent to
$a$, a contradiction to the assumption that $K$ is a module of
$G'$. Since $K-B' \not= \emptyset$, $K$ lies completely in $F_2$.
But now $B'$ is a join buoy of $G$, a contradiction.  $\Box$

\begin{Lemma}\label{lem:k5}
For an integer $h\ge 0$, let $H$ be any graph that is a pseudo-buoy with
bags $B_1, \ldots, B_5$ such that, for each $i\bmod 5$, $B_i$ is
${k_i}$-colorable, where $k_1, \ldots, k_5$ are integers that satisfy
$k_i+k_{i+1}\le h$ for each $i\bmod 5$.  Then:
\begin{itemize}
\item[(i)]
If $\sum_{i=1}^5 k_i \le 2h$, then $\chi(H)\le h$.
\item[(ii)]
If $\sum_{i=1}^5 k_i > 2h$ and each $B_i$ is $k_i$-chromatic, then
$\chi(H) > h$.
\end{itemize}
\end{Lemma}
\noindent{\it Proof.} Proof of (i).   We establish property (i) by
induction on $h$.  Suppose $\sum_{i=1}^5 k_i \le 2h$.  If $h=0$,
the property holds trivially.  Now assume that $h\ge 1$. Suppose
some three bags of $H$ are empty. Then $H$ is disconnected or is
the join of some $B_i$ and $B_{i+1}$ and so $H$ is $h$-colorable
by the lemma's hypothesis. Thus, some three bags of $H$ must be
non-empty. Say that a pair $(i,i+1)$ is \emph{tight} if
$k_i+k_{i+1}=h$.  Say that a pair $(j,j+2)$ is \emph{good}  if
$B_j$ and $B_{j+2}$ are both not empty, and every tight pair
contains an element from $\{j,j+2\}$.

Suppose all five pairs are tight. Then, we have
$5h=\sum_{i=1}^{5}(k_i+k_{i+1})=2\sum_{i=1}^{5}k_i \leq 4h$, which
is impossible. So, there is a non-tight pair. Suppose there is no
good pair. Among all non-tight pairs, choose the pair $(i,i+1)$
such that $k_i + k_{i+1}$ is minimized. We may assume without loss
of generality that $(1,2)$ are such a non-tight pair. If $B_3 =
\emptyset$, then the choice of the pair $(1,2)$ implies $B_1 =
\emptyset$. Now the sets $B_2, B_4, B_5$ are not empty, and so
$(2,4)$ is a good pair, a contradiction. Now we may assume $B_3
\not= \emptyset$, and by symmetry $B_5 \not= \emptyset$. But now
$(3,5)$ is a good pair, a contradiction.

So there is a good pair.  Up to relabelling, assume that $(1,3)$
is a good pair.  Pick any $k_1$-coloring of $B_1$, and let $S'$ be
a color class in that coloring.  Likewise, pick any $k_3$-coloring
of $B_3$, and let $S''$ be a color class in that coloring.  Let
$S=S'\cup S''$ and $H^*=H\setminus S$.  Thus $H^*$ is a
pseudo-buoy with bags $B_1'=B_1\setminus S',\; B_2'= B_2,\; B_3'=
B_3\setminus S'', \; B_4'=B_4, \; B_5'=B_5$.   Let $k_i' =
\chi(B_i')$ for $i=1,2,\ldots,5$. The fact that $(1,3)$ is a good
pair implies that $H^*$ satisfies the induction hypothesis for the
integer $h-1$.  In particular, since $\chi(B_i') \leq k_i -1$ for
$i= 1,3$, we have $\sum_{i=1}^5 k_i' \le 2(h-1)$. So $H^*$ is
$(h-1)$-colorable, and consequently $H$ is $h$-colorable since we
may use the $h$-th color on the vertices of~$S$.

Proof of (ii).  Suppose on the contrary that $\chi(H)\le h$.  So the
vertices of $H$ can be partitioned into $h$ stable sets $S_1, \ldots,
S_{h}$.  For each $i\in\{1, \ldots, 5\}$, let $\sigma_i$ be the number
of sets among $S_1, \ldots, S_{h}$ that have non-empty intersection
with $B_i$.  The definition of a pseudo-buoy implies that each stable
set $S_j$ has non-empty intersection with at most two of $B_1, \ldots,
B_5$.  It follows that $\sum_{i=1}^5 \sigma_i\le 2h$.  So
$\sum_{i=1}^5 \sigma_i < \sum_{i=1}^5 k_i$, which implies that there
is an integer $i\in\{1, \ldots, 5\}$ such that $\sigma_i< k_i$.  Hence
the non-empty intersections of $S_1, \ldots, S_{h}$ in $B_i$ form a
coloring of $B_i$ with strictly fewer colors than $k_i$, a
contradiction.  $\Box$

\medskip

Let $\mathcal{C}_k$ be the family of $k$-critical {\PP}-free
graphs. Clearly, $\mathcal{C}_1=\{K_1\}$ and
$\mathcal{C}_2=\{K_2\}$.  In general we have $K_k\in
\mathcal{C}_k$.  We are now in a position to prove the main
theorem of this section.

\begin{Theorem}\label{thm:construction}
For any $k\ge 2$, a graph $G$ is in $\mathcal{C}_k$ if and only if
it can be obtained by any of the following two constructions: \\
--- Construction 1: $G$ is the join of a member of $\mathcal{C}_{k_1}$
and a member of $\mathcal{C}_{k_2}$ for positive integers $k_1$ and
$k_2$ such that $k_1+k_2=k$.  \\
--- Construction 2: $G$ is a buoy with bags $B_1, \ldots, B_5$ such
that, for each $i\bmod 5$, the graph $B_i$ is a member of
$\mathcal{C}_{k_i}$, where $k_1, \ldots, k_5$ are positive integers
that satisfy $k_i+k_{i+1}\le k-1$ for each $i\bmod 5$ and
$\sum_{i=1}^5 k_i = 2k-1$.
\end{Theorem}

\noindent{\it Proof.} Let $k\ge 2$.

(I) We prove that any graph obtained by Construction~1 or~2 is in
$\mathcal{C}_k$.

First suppose that $G$ is obtained by Construction~1, i.e., $G$ is the
join of a member $G_1$ of $\mathcal{C}_{k_1}$ and a member $G_2$ of
$\mathcal{C}_{k_2}$ for positive integers $k_1$ and $k_2$ such that
$k_1+k_2=k$.  By Lemma~\ref{lem:join-critical}, $G$ is $k$-critical.

Now suppose that $G$ is obtained by Construction~2, with the same
notation as in the theorem.  By Lemma~\ref{lem:k5}~(ii) (with
$h=k-1$), we have $\chi(G)\ge k$.  Pick any $x\in B_1$.  We know that
$B_1 - x$ is $(k_1-1)$-colorable.  Moreover, we have
$(k_1-1)+\sum_{i=2}^5 k_i=2(k-1)$.  So the graph $G - x$, which is the
buoy with bags $B_1 - x, B_2, B_3, B_4, B_5$, satisfies the hypothesis
of Lemma~\ref{lem:k5}~(i) with $h=k-1$, and consequently it is
$(k-1)$-colorable.  The same holds for every vertex $x$ in $G$.  This
implies that $\chi(G)=k$ and $G$ is $k$-critical.

(II) Now we prove the converse part of the theorem.  Let $G=(V,E)$ be
a member of $\mathcal{C}_k$, i.e., $G$ is $(P_5, \overline{P}_5)$-free
and $k$-critical.  If $\overline{G}$ is not connected, then we can
apply Lemma~\ref{lem:join-critical}, and it follows that $G$ is
obtained by Construction~1.  Therefore assume that $\overline{G}$ is
connected.  Thus $G$ is not a clique and by Observation~\ref{obs:perfect} and
Lemma~\ref{lem:CHMW}, $G$
contains a $C_5$. By Lemma~\ref{lem:buoy}, $G$ contains a buoy $B$ which is  full or join.

Suppose $B$ is a join buoy of $G$.  Consider the set $A= V -B$.  Since
$A$ is complete to $B$, by Lemma~\ref{lem:join-critical}, both $A$ and
$B$ are critical, and so $G$ is obtained by Construction~2.

So, we may assume $B$ is a full buoy, that is, it contains all
vertices of $G$.  Let the bags of $B$ be $B_1, \ldots, B_5$.  Let
$k_i=\chi(B_i)$ for each $i\in\{1,\ldots, 5\}$.  For each $i$ we must
have $k_i+k_{i+1}\le k-1$, for otherwise $\chi(G[B_i\cup B_{i+1}])=k$,
which contradicts the fact that $G$ is $k$-critical.  Also we must
have $\sum_{i=1}^5 k_i \ge 2k-1$, for otherwise, by
Lemma~\ref{lem:k5}~(i), $G$ is $(k-1)$-colorable.  Each $B_i$ (with at
least two vertices) is a module of $G$, and so by
Lemma~\ref{lem:module-critical}, $B_i$ is $k_i$-critical.  Suppose
that $\sum_{i=1}^5 k_i \ge 2k$.  Pick any $x\in V$.  Then $G - x$ is a
pseudo-buoy that satisfies Lemma~\ref{lem:k5}~(ii), so $\chi(G - x)
\geq k$, which contradicts the fact that $G$ is $k$-critical.  So we
have $\sum_{i=1}^5 k_i = 2k-1$.  This shows that $G$ can be obtained
by Construction~2.  $\Box$

\begin{Theorem}\label{thm:finite}
For every $k$, $\mathcal{C}_k$ is a finite set.
\end{Theorem}
{\it Proof.} By induction on $k$. When $k=3$, it is easy to see
$\mathcal{C}_3$ contains two graphs: the $C_3$ and $C_5$ (every
graph $G$ in $\mathcal{C}_3$ must contain an odd cycle, and if $G$
does not contain $C_3$ or $C_5$, then $G$ contains an odd
chordless cycle with at least seven vertices and thus a $P_5$.)
Let $\mathcal{J}_k$ be the set of graphs of ${\cal C}_k$
constructed by Construction 1. Let $\mathcal{B}_k$ be the set of
graphs of $\mathcal{C}_k$ constructed by Construction 2. We have
$\mathcal{C}_k = \mathcal{J}_k \cup \mathcal{B}_k$ by
Theorem~\ref{thm:construction}.  Let $f(k)$ (respectively, $j(k)$,
$b(k)$) be the cardinality of $\mathcal{C}_k$ (respectively,
${\cal J}_k$, $\mathcal{B}_k$). Each graph in $\mathcal{J}_k$ is
constructed by taking the join of a graph in $\mathcal{C}_i$ and a
graph in ${\cal C}_{k-i}$ for $i = 1, 2, \ldots, k-1$.  It follows
that $j(k) \leq \sum_{i=1}^{k-1} f(i) f(k-i)$.  Consider a graph
$G$ in ${\cal B}_k$ which is a full buoy.  Graph $G$ has five
bags, each of which is a graph of $\mathcal{C}_i$ with $i \leq
k-2$.  It follows that $b(k) \leq (f(k-2))^5$.  Since $f(k) = j(k)
+ b(k)$, we have $f(k) \leq \sum_{i=1}^{k-1} f(i) f(k-i) +
(f(k-2))^5$.  So $f(k)$ is a function in $k$ and the result
follows.

For completeness, we will give a bound for $f(k)$ using Knuth's
up-arrow notation.  Define a single up-arrow operation to be $a
\uparrow b = a^b$.  Next define a double arrow operation to be $a
\uparrow \uparrow b = a \uparrow (a \uparrow (\ldots \uparrow a))$,
that is, a tower of $b$ copies of $a$.  It is easy to prove by
induction that $f(k) \leq 5 \uparrow \uparrow k$.
$\Box$
\subsection{A certifying algorithm for $k$-colorability of $(P_5,
\overline{P}_5)$-free graphs}

The algorithm in \cite{HoaLaz2013} accepts as input any graph $G$
that is $(P_5, \overline{P}_5)$-free and determines whether $G$ is
$k$-colorable in time $O(n^3)$, where $n$ is the number of
vertices of $G$.  However, when $G$ is not $k$-colorable, the
algorithm does not produce an easily verifiable certificate for a
``NO'' answer.

Now, given a $(P_5, \overline{P}_5)$-free graph $G$, we have a
simple polynomial-time algorithm for finding a $k$-critical
induced subgraph $H$ of $G$, as follows.  Consider a vertex $x$,
and determine if $G\setminus x$ is $k$-colorable; if it is not
$k$-colorable, then we remove $x$ from consideration; if it is
$k$-colorable, then $x$ is in~$H$.  Repeat the process for all
other vertices. The vertices that have not been removed from
consideration form the desired graph $H$. Now it is easy to check
whether $H$ satifies the properties given in
Theorem~\ref{thm:construction}, because, by
Theorem~\ref{thm:finite}, there is only a finite number of
$k$-critical graphs.  So it takes constant time to check whether
$H$ belongs to the family $\mathcal{C}_k$ of $k$-critical graphs.

\section{$k$-critical $(P_5, \overline{P}_5)$-free graphs for small
$k$}\label{sec:5-critical}

By Theorem~\ref{thm:finite} the number of $k$-critical {\PP}-free
graphs is finite for every $k$.  In this section, we will refine the
argument to establish sharp bounds on $|\mathcal{C}_k|$ for small
values of $k$.  In particular, we will construct all 5-critical
{\PP}-free graphs.

Let classes $\mathcal{B}_k$ and $\mathcal{J}_k$ be defined as in
Theorem~\ref{thm:finite}.  Recall that $\mathcal{C}_1=\{K_1\}$, ${\cal
C}_2=\{K_2\}$ and $\mathcal{C}_3=\{K_3, C_5\}$, so
$\mathcal{B}_1={\cal B}_2=\emptyset$ and $\mathcal{B}_3=\{C_5\}$.

When $\mathcal{A}$ and $\mathcal{B}$ are two sets of graphs, let
${\cal A}\otimes\mathcal{B}$ be the set of graphs that are the join of
a member of $\mathcal{A}$ and a member of $\mathcal{B}$.  We know that
$\mathcal{J}_k = \bigcup_{k_1+k_2=k} \mathcal{C}_{k_1}\otimes {\cal
C}_{k_2}$.  This means that each member $G$ of $\mathcal{J}_k$ is
either the join of several buoys or the join of $K_p$ and a member of
$\mathcal{C}_{k-p}$ for some positive $p$; in the latter case $G$ is
also the join of $K_1$ and a member of $\mathcal{C}_{k-1}$.  It
follows that we can write
\begin{equation}\label{eq:Jk}
\mathcal{J}_k = (\mathcal{C}_1\otimes \mathcal{C}_{k-1}) \cup (\bigcup
\mathcal{B}_{k_1}\otimes \cdots\otimes \mathcal{B}_{k_p})
\end{equation}
where the union is over all vectors $(k_1, \ldots, k_p)$ such that
$k_1+\cdots+k_p=k$ and $k_i\ge 3$.

When $G$ is a member of $\mathcal{B}_k$ we associate with $G$ its
\emph{pattern}, which is the numerical vector $(k_1, \ldots, k_5)$
(with the same notation as in Theorem~\ref{thm:construction}) that
satisfies the~constraints $k_i>0$, $k_i+k_{i+1}\le k-1$ for each
$i\bmod 5$ and $\sum_{i=1}^5 k_i = 2k-1$.  Several non-isomorphic
members of $\mathcal{B}_k$ can have the same pattern.  Conversely,
given a pattern $(k_1, \ldots, k_5)$ that satisfies the constraints,
one can construct a member of $\mathcal{B}_k$ as a buoy $(B_1, \ldots,
B_5)$ where $B_i$ is chosen from $\mathcal{C}_{k_i}$ for each $i$, and
Theorem~\ref{thm:construction} says that every member of
$\mathcal{B}_k$ is constructed that way (different choices may yield
the same member of $\mathcal{B}_k$ due to isomorphism).  We illustrate
this for the values $k=4$, $k=5$ and $k=6$.
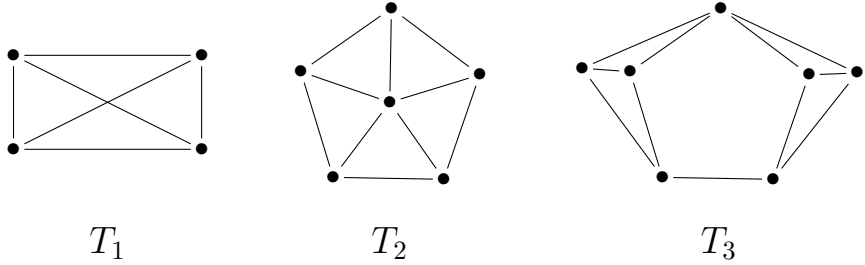
\begin{figure}
\begin{center}
\begin{tikzpicture} [scale = 1.25]
\tikzstyle{every node}=[font=\small]

\newcommand{\adj}{17}

\path (0,0) coordinate (g1);

\path (g1) +(-1,0.5) node (g1_1){}; \path (g1) +(1,0.5) node
(g1_2){}; \path (g1) +(-1,-0.5) node (g1_3){}; \path (g1)
+(1,-0.5) node (g1_4){};

\foreach \Point in {(g1_1),(g1_2),(g1_3),(g1_4)}{
    \node at \Point {\textbullet};
}

\draw   (g1_1) -- (g1_2)
        (g1_1) -- (g1_3)
        (g1_1) -- (g1_4)
        (g1_2) -- (g1_3)
        (g1_2) -- (g1_4)
        (g1_3) -- (g1_4);

\path (g1) ++(0,-1.5) node[draw=none,fill=none] { {\large $T_1$}};

\path (3,0) coordinate (g2);

\path (g2) +(0+\adj:1) node (g2_1){}; \path (g2) +(72+\adj:1) node
(g2_2){}; \path (g2) +(144+\adj:1) node (g2_3){}; \path (g2)
+(216+\adj:1) node (g2_4){}; \path (g2) +(288+\adj:1) node
(g2_5){}; \path (g2) node (g2_6){};

\foreach \Point in {(g2_1),(g2_2),(g2_3),(g2_4),(g2_5),(g2_6)}{
    \node at \Point {\textbullet};
}

\draw   (g2_1) -- (g2_2)
        (g2_2) -- (g2_3)
        (g2_3) -- (g2_4)
        (g2_4) -- (g2_5)
        (g2_5) -- (g2_1)
        (g2_6) -- (g2_1)
        (g2_6) -- (g2_2)
        (g2_6) -- (g2_3)
        (g2_6) -- (g2_4)
        (g2_6) -- (g2_5);
\path (g2) ++(0,-1.5) node[draw=none,fill=none] { {\large $T_2$}};

\path (6.5,0) coordinate (g3);

\path (g3) +(0+\adj:1) node (g3_1){}; \path (g3) +(72+\adj:1) node
(g3_2){}; \path (g3) +(144+\adj:1) node (g3_3){}; \path (g3)
+(216+\adj:1) node (g3_4){}; \path (g3) +(288+\adj:1) node
(g3_5){}; \path (g3) +(-5+\adj:1.5) node (g3_6){}; \path (g3)
+(149+\adj:1.5) node (g3_7){};

\foreach \Point in
{(g3_1),(g3_2),(g3_3),(g3_4),(g3_5),(g3_6),(g3_7)}{
    \node at \Point {\textbullet};
}

\draw   (g3_1) -- (g3_2)
        (g3_2) -- (g3_3)
        (g3_3) -- (g3_4)
        (g3_4) -- (g3_5)
        (g3_5) -- (g3_1)
        (g3_6) -- (g3_1)
        (g3_6) -- (g3_2)
        (g3_6) -- (g3_5)
        (g3_7) -- (g3_2)
        (g3_7) -- (g3_3)
        (g3_7) -- (g3_4);
\path (g3) ++(0,-1.5) node[draw=none,fill=none] {{\large $T_3$}};
\end{tikzpicture}

\end{center}
\caption{All 4-critical {\PP}-free graphs}\label{fig:4critical}
\end{figure}

For $k=4$, the only pattern that satisfies the constraints
$k_i>0$, $k_i+k_{i+1}\le 3$ for each $i\bmod 5$ and $\sum_{i=1}^5
k_i = 7$ is, up to circular permutation, $(1, 2, 1, 2, 1)$.  The
corresponding member of $\mathcal{B}_4$ is graph $T_3$ in
Figure~\ref{fig:4critical}.  By (\ref{eq:Jk}), we have
$\mathcal{J}_4 = \mathcal{C}_1\otimes \mathcal{C}_3$, so the
members of $\mathcal{J}_4$ are the graphs $T_1$ and $T_2$ in
Figure~\ref{fig:4critical}.  Hence $|\mathcal{C}_4|=3$.
\begin{figure}
\begin{center}
\input{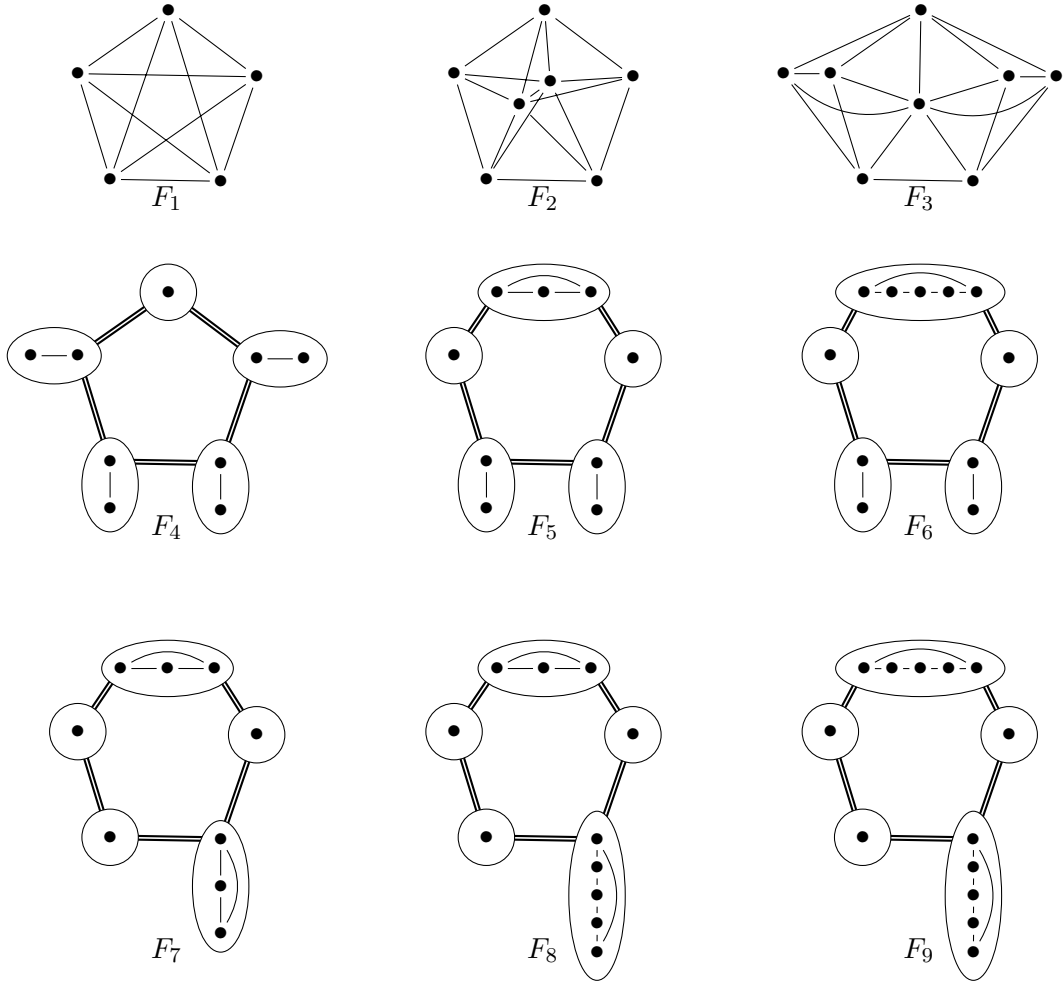}
\end{center}
\caption{All 5-critical {\PP}-free graphs.  In the graphs $F_4,
\ldots, F_9$, the ovals represent the bags and the double line denotes
all edges between the two bags.}\label{fig:list}
\end{figure}

For $k=5$, the patterns that satisfy the constraints $k_i>0$,
$k_i+k_{i+1}\le 4$ for each $i\bmod 5$ and $\sum_{i=1}^5 k_i = 9$ are,
up to circular permutation, $(2, 2, 2, 2, 1)$, $(2, 1, 3, 1, 2)$ and
$(1, 3, 1, 3, 1)$.  Pattern $(2,2,2,2,1)$ yields the graph $F_4$ on
Figure~\ref{fig:list}.  Since $\mathcal{C}_3=\{K_3, C_5\}$, pattern
$(2, 1, 3, 1, 2)$ yields graphs $F_5$ and $F_6$, and pattern
$(1,3,1,3,1)$ yields graphs $F_7$, $F_8$ and $F_9$.  So
$|\mathcal{B}_5|=6$.  By (\ref{eq:Jk}), we have $\mathcal{J}_5 = {\cal
C}_1\otimes \mathcal{C}_{4}$, so $\mathcal{J}_5$ consists of graphs
$F_1$, $F_2$, $F_3$ on Figure~\ref{fig:list}.  Hence $|{\cal C}_5|=9$.

For $k=6$, the patterns that satisfy the constraints $k_i>0$,
$k_i+k_{i+1}\le 5$ for each $i\bmod 5$ and $\sum_{i=1}^5 k_i = 11$
are, up to circular permutation, $(2, 2, 2, 2, 3)$, $(2,3,2,3,1)$,
$(2,3,1,3,2)$, $(1,4,1,3,2)$ and $(1,4,1,4,1)$.  Since $|{\cal
C}_1|=1$, $|\mathcal{C}_2|=1$, $|\mathcal{C}_3|=2$ and
$|\mathcal{C}_4|=3$, we see that $(2, 2, 2, 2, 3)$ yields two graphs,
$(2,3,2,3,1)$ yields four graphs, $(2,3,1,3,2)$, yields three graphs,
$(1,4,1,3,2)$ yields six graphs, and $(1,4,1,4,1)$ yields, up to
symmetry, six graphs.  So $|\mathcal{B}_6|=21$.  By (\ref{eq:Jk}), we
have $\mathcal{J}_6 = (\mathcal{C}_1\otimes \mathcal{C}_{5})\cup
({\cal B}_3\otimes \mathcal{B}_3)$, so $|\mathcal{J}_6|=9+1=10$.
Hence $|{\cal C}_6|=31$.

Similar computations lead to $|\mathcal{C}_7|=185$ and $|{\cal
C}_8|=1487$.


\begin{center}
{\bf Acknowledgement}
\end{center}
This work was done by authors H.S.D., T.J.D.M., and S.A.P. in partial
fulfillment of the course requirements for CP493: Directed Research
Project I in the Department of Physics and Computer Science at Wilfrid
Laurier University.  The authors A.M.H. and C.T.H. were each supported
by individual NSERC Discovery Grants.  The visit of author F.M. to
Wilfrid Laurier University was also supported by the NSERC Discovery
Grant of author C.T.H.

\end{document}